# DPNO: A Dual Path Architecture For Neural Operator

YiChen WANG[1], WenLian LU[1,2]

1 ShangHai Center of Mathematics, Fudan University, Shanghai, China
2 Department of Mathematics, Fudan University, Shanghai, China
23110840012@m.fudan.edu.cn, wenlian@fudan.edu.cn

**Abstract.** Neural operators have emerged as a powerful tool for solving partial differential equations (PDEs) and other complex scientific computing tasks. However, the performance of single operator block is often limited, thus often requiring composition of basic operator blocks to achieve better performance. The traditional way of composition is staking those blocks like feedforward neural networks, which may not be very economic considering parameter-efficiency tradeoff. In this paper, we propose a novel dual path architecture that significantly enhances the capabilities of basic neural operators. The basic operator block is organized in parallel two paths which are similar with ResNet and DenseNet. By introducing this parallel processing mechanism, our architecture shows a more powerful feature extraction and solution approximation ability compared with the original model. We demonstrate the effectiveness of our approach through extensive numerical experiments on a variety of PDE problems, including the Burgers' equation, Darcy Flow Equation and the 2d Navier-Stokes equation. The experimental results indicate that on certain standard test cases, our model achieves a relative improvement of over 30% compared to the basic model. We also apply this structure on two standard neural operators (DeepONet and FNO) selected from different paradigms, which suggests that the proposed architecture has excellent versatility and offering a promising direction for neural operator structure design.

## 1  Introduction

Partial differential equations (PDEs) are instrumental in modeling and analyzing a wide array of phenomena across physics, engineering and finance. In physics, wave propagation and heat diffusion are modeled by PDEs, whose solutions are essential for acoustics and thermal engineering. Besides, fluid dynamics also relies on PDEs like the Navier-Stokes equations to study turbulent flow which is crucial for engineering designs and meteorology. In Engineering, applications include structural mechanics for stress analysis and heat transfer in system design both rely on the effective solution of depended PDEs. Finance also employs PDEs in the Black-Scholes model for option pricing, impacting economic modeling. Thus, PDEs are a cornerstone in scientific and industrial advancements, offering solutions to complex, dynamic problems.

Due to the complexity of partial differential equation problems, it is difficult to obtain their analytical solutions in most cases; thus, it is necessary to rely on numerical methods to obtain numerical solutions which enables simulations, predictions, and system optimization. In practice, tasks such as industrial simulations, weather forecasting, geophysical exploration, and hydrological feature exploration place higher demands on the speed of solving partial differential equations. Therefore, the development of faster numerical methods as surrogate models has become an urgent problem to be addressed. However, traditional PDE solvers like finite difference (FDM) or finite element methods (FEM), due to its dependence on grid size, are computationally expensive, especially for high dimensional problems, and may struggle with complex geometries or multi-scale dynamics, limiting their flexibility and efficiency in modern applications.

With the widespread application of neural network methods in the scientific domain, numerous neural network-based methods for solving partial differential equations have been proposed. (such as [1, 2, 3, 4]) which addressed the limitations of traditional numerical methods. One prominent approach is the Deep ritz method proposed by [3], which is a powerful, mesh-free approach and is particularly effective in complex geometries and high dimensions. Another notable approach is Physics-Informed Neural Networks (PINNs) [4], which combines data-driven learning with physics-based modeling by incorporating the PDE residual into the loss function. Specifically, PINN learns the solution by minimizing the violation of the physical laws, i.e. the governing equation.

Although using neural networks to approximate a single partial differential equation have demonstrated superiority over traditional numerical methods to some extent, they still largely rely on grid partitioning. Additionally, different networks need to be trained for varying parameters, which significantly increases computational and application costs. Against this backdrop, a novel neural operator method has been proposed (such as [5], [6]). This approach focuses on learning mappings between function spaces. Once trained, it can handle different parameter configurations, and, through specially designed architectures (like [7]), can achieve resolution invariance, which means that it can be trained on a low resolution data and be applied to predict high resolution data.

Based on this foundation, our experiments reveal that a single-layer neural operator has limited representational capacity and could benefit from specially designed architectures. Indeed, in existing literature ([7, 8, 9]) neural operators are typically designed in a stackable form, and experimental results have demonstrated that this approach achieves superior performance compared to individual neural operators. Inspired by these results and recent advances in neural network architecture design within the field of image classification [10], we proposed a new dual path architecture that can be applied to basic neural operators. By comparing the newly designed model with the original model, we found that this new structure significantly improves the model's performance, with relative improvements exceeding 30% in some test cases.

## 2 Related Works

### 2.1 Neural Operator

Neural operators are a class of machine learning models designed to learn mappings between infinite-dimensional function spaces, making them particularly effective for solving partial differential equations (PDEs) and modeling complex physical systems. Unlike traditional neural networks that operate on finite-dimensional inputs, neural operators generalize across different resolutions and domains, enabling efficient predictions for unseen scenarios ([11]).

From the perspective of architectural design, neural operators can broadly be categorized into two fundamental paradigms: the encoder-decoder paradigm and the kernel integral operator paradigm. The first paradigm, based on the universal approximation theorem of nonlinear operators [12], employs neural networks to approximate the basis of the target function space, which then serves as the approximation for the final function. In contrast, the second paradigm, proposed by [5], leverages the concept of Green's functions, treating the target operator as a composition of simple integral kernel operators. The first paradigm typically exhibits a more independent structure, while the second paradigm inherently possesses a layered structure. Since they were proposed, both of these paradigms have each spawned a significant amount of related research ([2, 5, 6, 7, 13, 14, 15]). To demonstrate the general applicability of our method, we select one representative model from each paradigm and compared them with their redesigned version.

### 2.2 Network Structure

In the study of neural networks, the design of neural network architectures has always been an important topic and have evolved significantly since the early 2000, driven by research aimed at improving performance and efficiency. Several well-known examples include U-Net [16], ResNet [17], DenseNet [18], and ViT [19].

U-Net is a deep learning architecture designed primarily for image segmentation tasks. Introduced by Ronneberger et al. in 2015 ([16]), the U-Net structure consists of an encoder (contracting path) and a decoder (expansive path), connected by a bottleneck. The encoder downsamples the input image through convolutional layers, capturing context, while the decoder upsamples the feature maps to reconstruct the original resolution. Skip connections bridge the encoder and decoder, merging low-level features from the encoder with high-level features from the decoder, enabling precise localization and accurate segmentation. This architecture balances context and spatial information, making it highly efficient for segmentation tasks.

Residual Networks (ResNet) introduced residual connections, enabling the training of deeper networks by mitigating vanishing gradients. This breakthrough allows for more complex models, enhancing accuracy in tasks like image classification. Building on this, DenseNet further optimized feature reuse by connecting each layer to every other, reducing parameter redundancy and improving efficiency.

Dual Path Networks (DPN), introduced in 2017 ([10]), proposed a dual path architecture that bridges the gap between two popular network designs: ResNet and DenseNet. By viewing them as higher order revolutionary neural networks, DPN employs two parallel paths for maintaining consistent feature learning and exploring new feature representations. Its efficiency in both feature learning and computational resource utilization makes it a compelling choice for modern deep learning applications.

Vision Transformers (ViT) ([19]) represents a distinctly different architectural paradigm from neural networks. By leveraging self-attention mechanisms to capture long-range dependencies, ViT demonstrates superior performance in various vision tasks, which provides new insights into architectural design.

### 2.3 Operator Structure

Building on numerous studies in network architecture design, recent research has proposed a variety of novel neural operators, which have achieved remarkable results in standard tasks; many of which attained state-of-the-art (SOTA) performance at the time of their introduction.

In literature [8], the author proposed an integral autoencoders for discretization-invariant learning, in which multi-channel learning was used to improve the performance of basic integral autoencoders. Specifically, they use a densely connected multi-block structure inspired by DenseNet. Their experiments show that this structure achieves a start-of-art performance. Enlightened by their work, we further explore the construction of multi-block structure and find a dual path architecture that can enhance the performance of basic operator even further.

Another work focus on the design of architecture of neural operator is the UNO [20], which is inspired by the U-Net. This work designed a U-shaped memory enhanced architecture that allows for deeper neural operators. UNO exploits the problem structures in function predictions and demonstrate fast training, data efficiency, and robustness with respect to hyperparameters choices.

Due to the scalability of the Transformer architecture and its superior performance demonstrated in various tasks, a significant body of research has begun to design neural operators based on the Transformer architecture. Notable works include: Galerkin Transformer [21], GNOT [22], ONO [23] and Transolver [9]. Galerkin Transformer is the first research that apply self-attention mechanism to data-driven operator learning problem, which proposes a new layer normalization scheme called Galerkin-type attention. After that, in 2023, GNOT was proposed focusing on the resolution of irregular grids, multi-input functions, and multi-scale problems in operator learning tasks. Through the designing of a Heterogeneous Normalized Attention block and a soft gating mechanism based on mixture-of-experts, GNOT effectively addresses these issues and significantly surpasses previous approaches in operator learning. In 2024, based on the idea of using Transformer to do kernel operator learning, the researchers constructed ONO and Transolver structure. The distinction lies in the fact that the former introduces an orthogonal attention mechanism, employing two non-coupled pathways to respectively approximate eigenfunctions and the solution of the equation. In contrast, the latter focuses on utilizing an attention mechanism to extract intrinsic physical states from discrete geometric information. Given that their design philosophy also involves achieving better results through the stacking of basic operator modules, therefore, the dual path architecture we propose is also applicable to these works.

## 3 Method

### 3.1 Problem Setup

In methodology, we consider PDEs as a mapping between two infinite dimensional spaces. Let $D \subset \mathbb{R}^d$ be a bounded open set in Euclidean space, then we denote the input function input function by $a(x)$, which can be the boundary condition, initial condition or parameter function appeared in the equation. The corresponding function space is denoted by $a(x) \in \mathcal{A}(D, R^{d_a})$. The output function is denoted by

$u(y)$, which is generally the solution of the PDE, and the output function space is $u(y) \in \mathcal{U}(D, R^{d_u})$. Let $\mathcal{G}^\dagger: \mathcal{A} \to \mathcal{U}$ be the underlying solution operator, i.e. $u(y) = \mathcal{G}^\dagger(a(x))(y)$. Our goal is to approximate $\mathcal{G}^\dagger$ by building a parametric neural network model $\mathcal{G}(\cdot, \theta)$, where $\theta \in \Theta$ is the learnable parameter. Then, we also need to define a cost functional $C: \mathcal{U} \times \mathcal{U} \to R$ to transform this problem into a network learning problem. After that the best approximator of the solution operator is given by $\mathcal{G}(\cdot, \theta^*)$, where $\theta^*$ is the minimizer of the following probelm:

$$\min_{\theta \in \Theta} \mathbb{E}_{a \sim \mu}[C(\mathcal{G}(a, \theta), \mathcal{G}^\dagger(a)]$$

In practice, we sample input function $a(x)$ independently from a Gaussian random field and solve the partial differential equations numerically to obtain the corresponding true solution $u(y) = \mathcal{G}^\dagger(a(x))(y)$, together, we have observation $\{a_j, u_j\}_{j=1}^N$. Further, since both $a_j$ and $u_j$ are functions, in order to deal them numerically, discretization is necessary. Let $D_j = \{x_1 \cdots, x_n\} \subset D$ be a n-points discretization, then, our input-output pairs are $\{a_j|_{D_j}, u_j|_{D_j}\}_{j=1}^N$.

### 3.2 Dual Path Architecture

Similar with dual path networks [10], the general dual path architecture is made of two distinct paths which are designed to capture more complex operator structure in the transition process. As shown in Fig. 1.

The first path is formed in a manner similar to ResNet. The equation is given by:

$$u_{k+1}(x) = \mathcal{G}_k(u_k)(x) + u_k(x)$$

The second path, inspired by DenseNet, uses densely connects operator blocks, allowing functions to be shared and reused across different depths. The equation is given by:

$$v_{k+1}(x) = \mathcal{G}_k([v_0, v_1, \cdots, v_k])(x)$$

Finally, the output of both paths are fed into a neural network to generate the final solution $u(x)$.

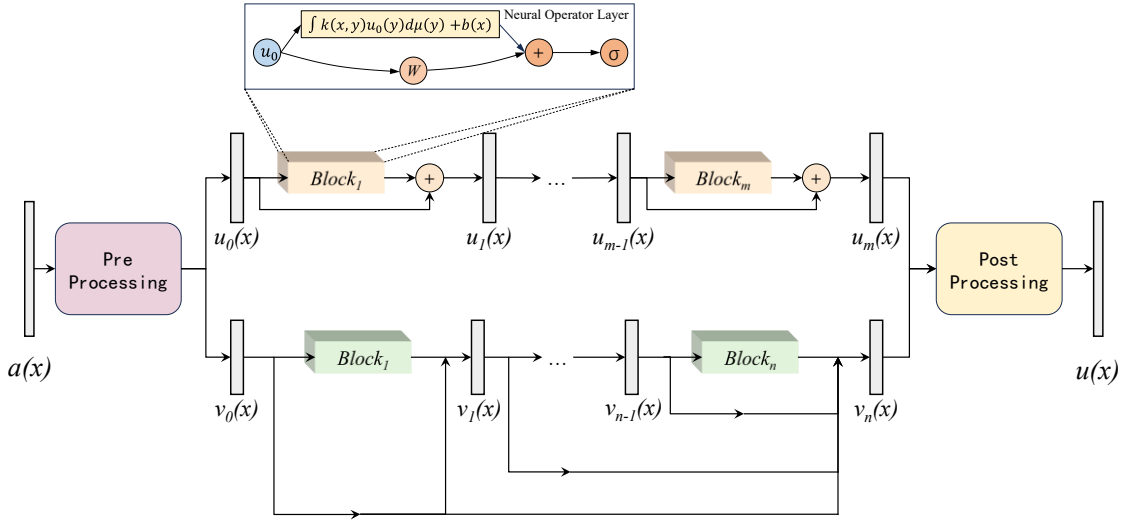

**Fig. 1.** Dual Path Architecture. This figure shows the dual path architecture, where $a(x)$ is the input function, $u(x)$ is the output function, $u_i(x), v_i(x)$ denotes the transition functions and blocks denote the operator block, which can be either Neural Operator or DeepONet, the figure displays Neural Operator as an example.

### 3.3 Operator Blocks

Although the problem formulations are identical, as mentioned above, two distinct paradigms have emerged on how neural networks can approximate map between infinite-dimensional function spaces ([24]): one based on Encoder-Decoder ([2, 6, 14]) and the other views operator as the composition of kernel integral operator ([5, 7, 9]).

To demonstrate the superiority of the dual path architecture, from each of the two paradigmatic categories, we selected one canonical model (DeepONet and FNO) as a representative and compared it with its dual path version. Due to differences in model construction, the dual path architecture can't be uniformly applied to two paradigms and must be specially designed according to the construction of the base models. Below, we will elaborate on the specific definitions of each.

### 3.4 Dual Path DeepONet

DeepONet was first proposed by Lu Lu et al. in [6], whose structure is designed based on the universal approximation theorem of nonlinear operators [12].

A basic DeepONet consists of two networks (named branch net and trunk net) that encode input function and grids for output function respectively.

The equation is given by (according to [6]):

$$\begin{cases} branch(a(x)) = \sum_{i=1}^{n} \sigma\left(\sum_{j=1}^{m} \xi_{ij} a(x_j) + \theta_i\right) \\ trunk(y) = \sigma(w \cdot y + \zeta) \\ \mathcal{G}(a)(y) \approx branch(a(x)) \cdot trunk(y) = \sum_{k=1}^{p} \sum_{i=1}^{n} c_i^k \sigma\left(\sum_{j=1}^{m} \xi_{ij}^k a(x_j) + \theta_i^k\right) \cdot \sigma(w_k \cdot y + \zeta_k) \end{cases}$$

where $a$ denotes the input function, $x_j$, $y$ denote the input and output grid, $\sigma$ denotes activation function, $\xi_{ij}^k, \theta_i^k$ are the parameters of branch net, $w_k, \zeta_k$ denotes the parameters of trunk net, $c_i^k$ is the summation coefficient when combining branch net and trunk net.

In the framework of the universal approximation theory, the output of trunk net can be viewed as bases function of output function space while the output of branch net can be viewed as the expansion coefficients of solution u(y) with respect to those bases. Therefore, in order to utilize the dual path structure and decrease the number of model simultaneously, we only apply dual path architecture on trunk net.

In details, the trunk net is taken as a basic block and in ResNet path, the input of each block is added to the output to formulate the input of subsequent block. In DenseNet path, the input of each block is concatenate with the output. This can be explained as increasing the number of basis functions used in approximating the solution. The experimental result validates that this explanation and shows a significantly enhanced performance compared with the original model.

### 3.5 Dual Path Fourier Neural Operator

FNO is an example of general operator architecture designed by Li et al. in [7], to describe its structure, we need some definitions:

**Definition 1: Fourier Transform**

The Fourier Transform $\mathcal{F}$ and its inverse $\mathcal{F}^{-1}$ are given by:

$$(\mathcal{F}f)(k) = \int_D f(x) e^{-2i\pi\langle x,k\rangle} dx,$$

$$(\mathcal{F}^{-1}f)(k) = \int_D f(x)e^{2i\pi\langle x,k\rangle}\,dk$$

**Definition 2: Fourier integral operator**

The Fourier integral operator $\mathcal{K}$ is the defined by:

$$(\mathcal{K}(\phi)a_t)(x) = \mathcal{F}^{-1}\left(R_\phi \cdot \mathcal{F}(a_t)\right)(x) \qquad \forall x \in D$$

where $R_\phi$ is the Fourier transform of the kernel function $\kappa: D \to R^{d_a \times d_a}$ parameterized by $\phi$.

Then, we adopt the description from [7] to outline the process as follows: Firstly, the input function $a(x)$ and its grid are preprocessed by a lifting map $P$, which is usually a fully connected neural networks. Next, the function are updated iteratively following the below equation:

$$a_{t+1}(x) := \sigma\big(Wa_t(x) + (\mathcal{K}(\phi)a_t)(x)\big)$$

where $W: R^{d_a} \to R^{d_a}$ is a linear transformation and $\mathcal{K}$ is the Fourier integral operator defined above. After several update steps, the output $a_T(x)$ is projected to $u(x)$ by a local transformation $Q: R^{d_a} \to R^{d_u}$.

The whole process is described by:

$$\begin{aligned}
a_0(x) &:= P\big(x, a(x)\big) \\
\mathcal{K}(\phi)(a_t)(x) &:= \mathcal{F}^{-1}(R_\phi \cdot (\mathcal{F}(a_t))(x) \\
a_{t+1}(x) &:= \sigma\big(Wa_t(x) + (\mathcal{K}(\phi)a_t)(x)\big) \\
u(x) &= Q\big(a_T(x)\big)
\end{aligned}$$

where $a(x)$ is the input function, $P$ is the lifting map, $\mathcal{F}$ and $\mathcal{F}^{-1}$ denote Fourier transform and inverse Fourier transform respectively, $T$ denotes the number of iteration steps, $t = 0,\cdots,T-1$ denotes the index number of iteration and $Q$ is the local linear transformation.

In designing dual path structure for Fourier Neural Operator, we follow the mode in original work [7], in which a kernel integral operator together with a nonlinear activation function is taken as a basic block. Then plus operation and concatenate operation are used between basic blocks to formulate ResNet path and DenseNet path.

According to the theory proposed by [10], the DenseNet path generates more effective operator, meanwhile, the ResNet path better leverages the operator generated by lower layers, which are then combined to formulate a more powerful general neural operator.

## 4    Experiments

To demonstrate the effectiveness of our dual path architecture, we select three operator approximation problems (Burgers' Equation, Darcy Flow Equation and 2-d Navier-Stokes Equation) according to [7] and compared the performance of DeepONet, FNO, and their dual path version on these problems.

The settings of them are as follows: For DeepONet, we use a 4 layers block with width 128 in the middle and GeLU activation to construct both branch net and trunk net. For FNO, we stack four Fourier integral operator blocks and also use GeLU as activation function. For dual path DeepONet, we use the same setting as DeepONet as basic block. For dual path FNO, we use a single Fourier integral operator block as basic block.

Since the deepening of DenseNet path will significantly increases the number of parameters, considering the GPU memory constraints, we use 3 operator blocks in the DenseNet path and 4 operator blocks in ResNet path.

All models are trained with 1000 training samples and tested with 200 test samples. We train 1000 epochs and test every 20 epochs using AdamW optimizer starting with learning rate 0.001. The Experiments are conducted on a Nvidia 4090 GPU with 16GB memory. The results are recorded in Table 1

To visually demonstrate the differences in model performance, we randomly select a test instance, and display the absolute error of model prediction (Two cases are shown). For FNO, we use Darcy flow data (displayed in Fig. 4), and for DeepONet we use 2d Navier-Stokes data (displayed in Fig. 6).

| Model | Burgers | Darcy Flow | Navier-Stokes |
|---|---|---|---|
| DeepONet | 0.02267 | 0.07368 | 0.05824 |
| Dual Path DeepONet | 0.01369 | 0.06311 | 0.04933 |
| **Relative Promotion** | **39.58%** | **14.34%** | **15.30%** |
| FNO | 0.001466 | 0.04768 | 0.002594 |
| Dual Path FNO | 0.0009712 | 0.03017 | 0.002402 |
| **Relative Promotion** | **33.77%** | **36.73%** | **7.38%** |

**Table 1**: Performance Comparison. The relative L2 loss between true solution and model prediction is recorded.

## 4.1 Burgers Equation

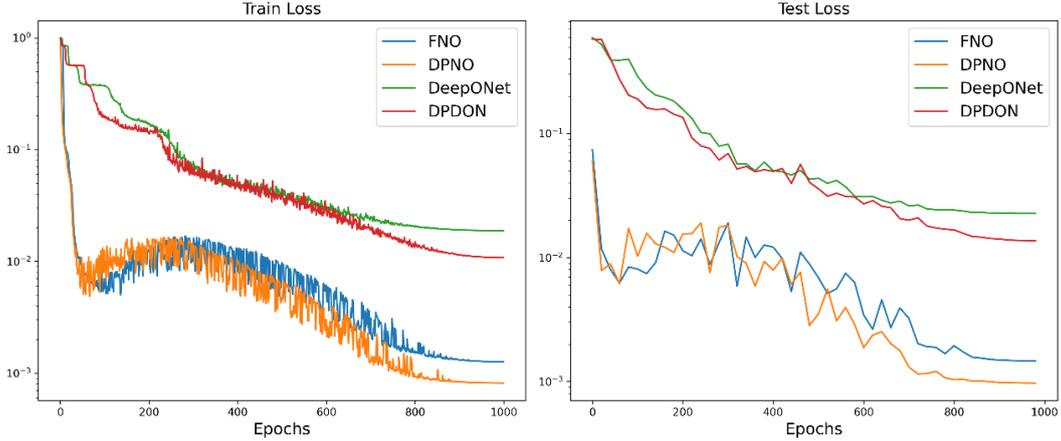

**Fig. 2.** Training and testing process of Burgers' Equation, x label is the epoch and y label is the relative $L^2$ loss. The test is conducted every 20 epochs.

For 1 dimensional case, we use Burger's Equation, which is a nonlinear equation originating from applying fundamental conservation laws (See [25] for details). This equation is of great importance in studying shock wave formation and turbulence, as well as in traffic flow modeling and biological systems.

In our experiment, we use the Burgers' equation on unit torus:

$$\begin{cases} \partial_t u(x,t) + \partial_x \left( \dfrac{u^2(x,t)}{2} \right) = \nu \partial_{xx} u(x,t), & t \in (0,1], \\ u(x,0) = u_0(x), & x \in (0,1). \end{cases}$$

where $u_0 \in L^2\big((0,1);\mathbb{R}\big)$ is the periodic initial condition generated according to $u_0 \sim \mathcal{N}(0, 625(-\Delta + 25I)^{-2})$ and viscosity $\nu$ is set by $\nu = 0.01$. The equation is solved on a spatial mesh with resolution 2048 and the dataset is sampled from this using a 256 resolution.

Our goal is learning the map from the initial condition $u_0(x)$ to the solution at time [0,1], i.e. $\mathcal{G}^\dagger: u_0(x) \to u(x,[0,1])$

## 4.2 Darcy Flow Equation

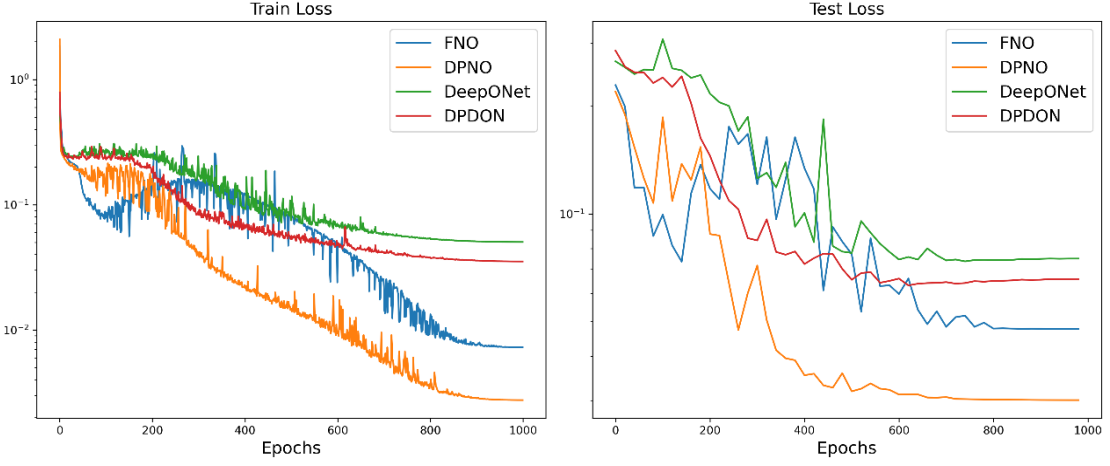

**Fig. 3.** Training and testing process of Darcy Flow Equation, x label is the epoch and y label is the relative $L^2$ loss. The test is conducted every 20 epochs.

For 2 dimensions instance, we consider Darcy Flow Equation on unit box with Dirichlet boundary, given by:

$$\begin{cases} \nabla \cdot (a(x)\nabla u(x)) = f(x), & x \in (0,1)^2 \\ u(x) = 0, & x \in \partial(0,1)^2. \end{cases}$$

where a(x) is the diffusion coefficient and f(x) is the forcing function.

Follow the setting in [7], in this experiment, $a \in L^\infty((0,1); R)$ is generated according to $a \sim \psi\left(\mathcal{N}\left(0, (-\Delta + 9I)^{-2}\right)\right)$ with zero Neumann boundary conditions on the Laplacian. The map $\psi(x): \mathbb{R} \to \mathbb{R}$ takes value 12 if $x \geq 0$, takes value 3 if $x < 0$. The forcing function $f(x)$ is taken to be $f(x) \equiv 1$.

The equation is solved on a spatial mesh with resolution 256*256 and the dataset is sampled from this using a 128*128 resolution. We aim to learn the operator mapping diffusion coefficient function a(x) to solution $u(x)$.

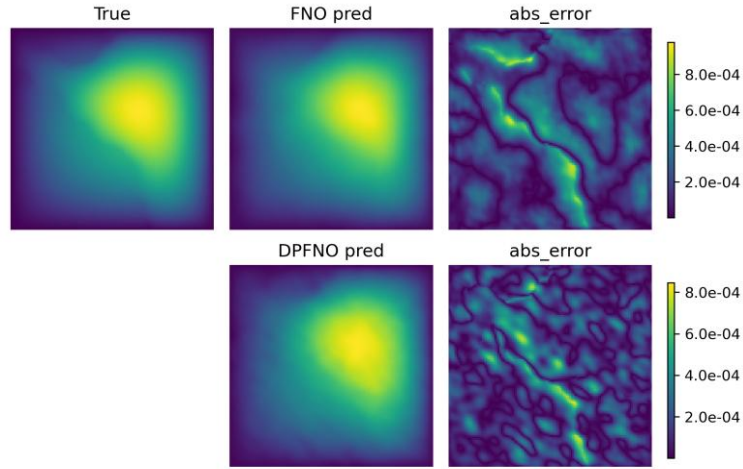

**Fig. 4.** Darcy Flow. This figure shows an instance of Darcy Flow equation. The first row, from left to right, are the true solution, the prediction made by FNO and its absolute error with respect to the true solution. The second row, are dual path FNO prediction and its absolute error with respect to the true solution.

## 4.3 Navier-Stokes Equation

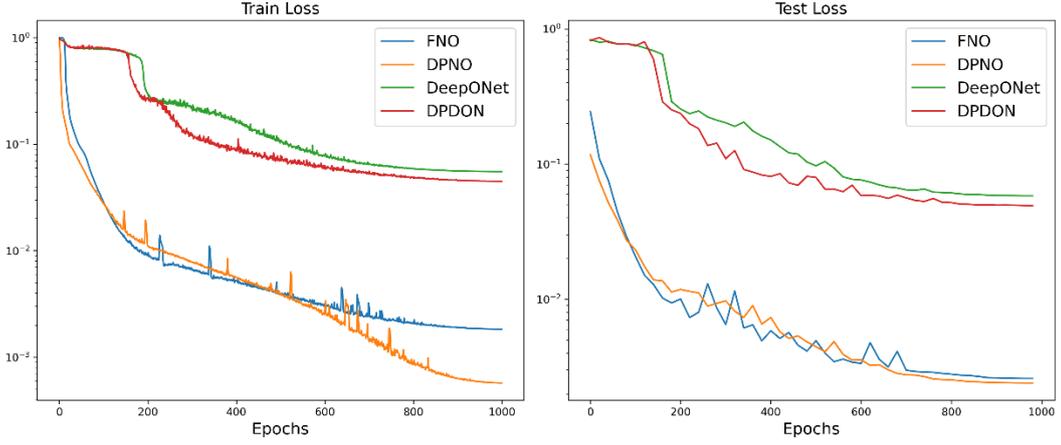

**Fig. 5.** Training and testing process of 2d Navier-Stokes Equation, x label is the epoch and y label is the relative $L^2$ loss. The test is conducted every 20 epochs.

For three dimensions problem, we consider Navier Stokes Equation, which is of immense importance since it provides the mathematical foundation for understanding phenomena like ocean currents and atmospheric flow.

In our experiment, we choose 2-d Navier Stokes Equation that models viscous, incompressible fluid in vorticity form on the unit torus.

The equation is given as follows:

$$\begin{cases} \partial_t w(x,t) + u(x,t) \cdot \nabla w(x,t) = \nu \Delta w(x,t) + f(x), \\ \nabla u(x,t) = 0, \quad x \in (0,1)^2, t \in (0,T] \\ w(x,t) = w_0(x), x \in (0,1)^2 \end{cases}$$

where $u$ is the velocity field, $w = \nabla \times u$ is the vorticity, $w_0 \in L^2((0,1)^2; R)$ is the initial vortcity generated according to $\mathcal{N}\left(0, 7^{3/2}(-\Delta + 49I)^{-2.5}\right)$, $\nu \in R_+$ is the viscosity coefficient, which is kept fixed 0.001 in our experiment. $f \in L^2$ is the forcing function, which is set to $0.1\left(sin(2\pi(x+y)) + cos(2\pi(x+y))\right)$. More details of dataset generation and can be found from [7].

We are interested in learning the operator between the vorticity at the first ten time steps to the vorticity at the next ten time steps, i.e. $\mathcal{G}^\dagger: w(x,(0,10]) \to w(x,[10,20])$.

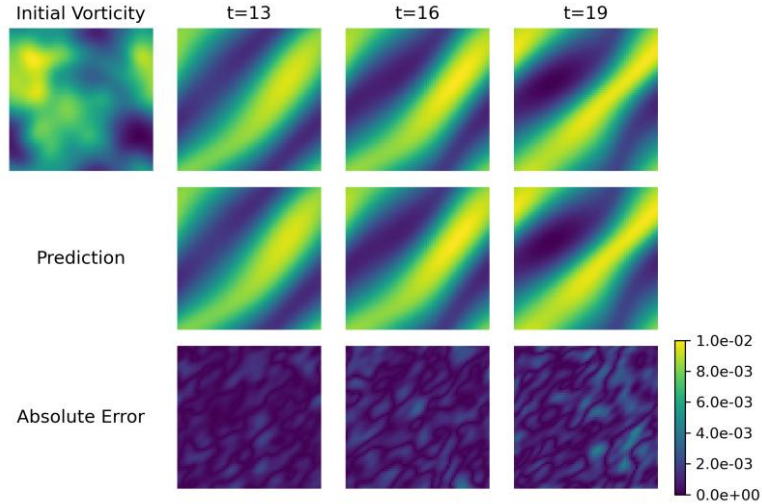

**Fig. 6.** Navier-Stokes. This figure shows the prediction made by Dual Path DeepONet on an instance of Navier-Stokes equation data. Due to space limitation, we only display 3 time steps. The first row, from left to right, are the initial vorticity, true solution at time 13, 16 and 19. The second row displays the corresponding prediction made by Dual Path DeepONet and the third row illustrates its absolute error with respect to the true solution.

## 5  Conclusions

This paper presents a dual path architecture that can be applied to basic neural operator blocks, which effectively compensates for the limitation of the single-block neural operator's representation ability. The experiments show an impressive performance improvement of dual path model compared with original model. Furthermore, this architecture can be applied to various fundamental neural operators, which demonstrates its excellent generalization capability. However, despite this model structure performs well in experiments, its effectiveness still lacks rigorous theoretical proof. Another noteworthy point is that we only modified two basic neural operator models (DeepONet and FNO). With the recent emergence of transformer-based neural operator models, the dual path architecture can also be applied as a design pattern to these new models. Further experiments are needed to supplement future work.